\documentstyle[11pt]{article}
\newmathalphabet*{\msbt}{cmr}{m}{sl}
\addtoversion{normal}{\msbt}{msb}{m}{n}
\addtoversion{bold}{\msbt}{msb}{m}{n}
\newmathalphabet*{\msat}{cmr}{m}{sl}
\addtoversion{normal}{\msat}{msa}{m}{n}
\addtoversion{bold}{\msat}{msa}{m}{n}

           \let\sec=\S    
           \let\norwegian=\O
	\let\pist=\P

 \def\a{\alpha}               \def\th{\theta}              \def\p{\pi}
 \def\b{\beta}                                
 \def\g{\gamma}               \def\k{\kappa}               
 \def\d{\delta}               \def\l{\lambda}              \def\s{\sigma} 
 \def\e{\epsilon}             \def\m{\mu}                  \def\t{\tau}
        \def\n{\nu}

                                                           \def\o{\omega}

                                \def\S{\Sigma}
 \def\D{\Delta}                  
                            \def\P{\Pi}
                                 \def\O{\Omega}

           \newcommand{\an}{\wedge}
                         \def\all{\forall}
           \def\Equi{\Longleftrightarrow}   \def\ex{\exists}
           \def\Imp{\Rightarrow}            
           \def\sset{\subseteq}		   
           \def\equi{\longleftrightarrow}
           
            \def\re{\mbox{$\msbt R$}}

           \def\int{{\msbt Z}}
           
          \newtheorem{thom}{Theorem}[section]
	  \newtheorem{fact}[thom]{Fact}
          \newtheorem{lemma}[thom]{Lemma}
          \newtheorem{propo}[thom]{Proposition}
          \newtheorem{de}[thom]{Definition}
          \newtheorem{corr}[thom]{Corollary}

	  \def\thm{\begin{thom}}           \def\ethm{\end{thom}}
          \def\lem{\begin{lemma}}	   \def\elem{\end{lemma}}		

         \def\prop{\begin{propo}}         \def\eprop{\end{propo}} 
          \def\df{\begin{de}}  		   \def\edf{\end{de}}
          \def\cor{\begin{corr}}           \def\ecor{\end{corr}}
	  \def\fa{\begin{fact}}  	   \def\ef{\end{fact}}       
          \def\bq{\begin{quote}}           \def\eq{\end{quote}}

\begin{document}
\def\text{\mbox}
\def\bm{\boldmath}
\def\pf{{\bf {Proof: }}}

\newcommand{\dl}{\mbox{\bm ${\d^1_2}$}}
\newcommand{\card}[1]{\mbox{$\mathrm{card}(#1)$}}
\newcommand{\qd}[1]{\hspace*{\fill{\bf \mbox{ Q.E.D.#1}}}}
\newcommand{\q}{\hspace*{\fill{\bf \mbox{ Q.E.D.}}}}
\newcommand{\rem}{\noindent {\bf{Remark: }}}
\newcommand{\Ka}{\mbox{$K^a$}}
\newcommand{\Um}{\mbox{${\cal U}$}}
\newcommand{\as}{a^\#}
\newcommand{\Kas}{\mbox{$K^{\as}$}}
\newcommand{\Las}{\mbox{$L[\as]$}}
\newcommand{\Kass}{\mbox{$K^{a^{\#\#}}$}}
\newcommand{\ul}{\mbox{$\underline{L}$}}
\newcommand{\uk}{\mbox{$\underline{K}$}}
\newcommand{\La}{\mbox{$L[a]$}}
\newcommand{\game}{G_{M,\a}}
\newcommand{\empseq}{\langle\,\rangle}
\newcommand{\rkempty}{\mbox{$rk_{T_{M,\a}}(\empseq)$}}

\newcommand{\omo}{\mbox{$\o_1$}}           
\newcommand{\omt}{\mbox{$\o_2$}}           
\newcommand{\om}[1]{\mbox{$\o_#1$}}       
\newcommand{\ba}{{\bar \a}}
\newcommand{\ta}{{\tilde \a}}
\newcommand{\tm}{\widetilde M}
\newcommand{\vd}{{\vec \d}}
\newcommand{\vg}{{\vec \g}}
\newcommand{\mi}{\bar M}
\newcommand{\Ni}{\bar N}
\newcommand{\pistol}{\scriptsize \pist}
\newcommand{\Nbar}{{\bar N}}
\newcommand{\Hbar}{{\bar H}}
\newcommand{\Kbar}{{\bar K}}
\newcommand{\thbar}{{\bar \th}}
\newcommand{\kbar}{{\bar \k}}
\newcommand{\Mbar}{{\bar M}}
\newcommand{\xbar}{{\bar x}}
\newcommand{\sbar}{{\bar \s}}
\newcommand{\pbar}{{\bar \p}}
\newcommand{\nbar}{{\bar \n}}
\def\K{$K$}
\def\M{$M$}
\def\N{$N$}
\def\W{$W$}
\def\Q{$Q$}
\def\wt{\widetilde}
\def\emp{\mbox{ \norwegian}}
\def\hright{\hookrightarrow}

\def\calt{{\cal T}}
\def\calv{{\cal V}}
\def\calw{{\cal W}}
\def\calu{{\cal U}}
\def\calj{{\cal J}}
\def\ite{{\calt}}

\hsize=16.5cm
\def\ptwo{\P^1_2}
\def\sit{\S^1_3}
\def\dit{\mbox{\bm $\d^1_2$}}
\newcommand{\fin}[1]{ [#1]^{<\o} }
\newcommand{\finn}[2]{ [#1]^{#2} }
\def\back{\backslash}
\def\rest{\restriction}
\def\la{\langle}
\def\ra{\rangle}
\def\dlim{\mbox{\raisebox{-1.25ex}{$\stackrel{\textstyle Lim}{\longrightarrow}\,$}}}
\input{amssym.def}
\undefine\upharpoonright
\newsymbol\upharpoonright 1316
  \newcommand{\restriction}{\mbox{$\upharpoonright$}}
\newsymbol\smallfrown 1361
\newsymbol\nleq 2302
\def\mar{\marginpar{}}

\newcommand{\mb}{\mbox}
\newcommand{\disp}{\displaystyle}
\newcommand{\beq}{\begin{equation}}
\newcommand{\eeq}{\end{equation}}
\newcommand{\und}{\underline}
\def\baire{\mbox{}^\o\o}
\def\dew{\D^1_3}
\def\wdew{\widetilde\D^1_3}
\def\wdewt{\widetilde\D^1_3(T)}
\def\dite{\cite}
\parindent=0pt
\parskip=.4cm

\thispagestyle{empty}

\setcounter{page}{1}
\title{
 Determinacy and $\mb{$\D$}^1_3$ - degrees}

\author{ P.D.Welch
\\School of Mathematics, University of Bristol.}
\date{July 16th 1996}

\maketitle

\begin{abstract}

Let $D =\{${\bm $d_n$}$\}$ be a countable collection of $\D^1_3$ degrees.
Assuming that all co-analytic games on integers are determined (or
equivalently that all reals have ``sharps''), we prove that either
$D$ has  a $\D^1_3$-minimal upper bound, or that for any $n$, and
for every real $r$ recursive in $d_n$, games in the pointclasses
$\D^1_2(r)$ are determined. This is proven 
using Core Model theory.

\end{abstract}

\newpage
\setcounter{page}{1}

\vspace{.2in}

 \section{Introduction}

The main theorem of \cite{W1} is the following:

\thm  $(ZFC + \forall r \in \re\;(r^{\#}$ exists$) + \neg \, 0^{\dagger}$)
  Every countable set of $\D^1_3$ degrees has a minimal upper
bound.
 \ethm 
We say that for reals $f,g \in \mbox{}^\omega\omega\ f$
{\em is} $\D^1_3$ {\em in} $g$, or $f \leq_3 g$, if there are
$\Sigma^1_3$ relations $\Phi, \Psi$ expressible in second order number
theory so that $f(n) = m \Leftrightarrow \Phi (n, m, g) \Leftrightarrow
\neg \Psi (n, m, g)$.  This is a reducibility ordering, and setting
$f = _3 g \Leftrightarrow f \leq_3 g\ \&\ g \leq_3 f$ we have that
$=_3$ is an equivalence on $\mbox{}^\omega\omega\ $.  We let
$\mb{\bm$f$}$ denote the equivalence class $[f]_{= 3}$.

In general, when $\leq_r$ is a reducibility ordering (meaning $\leq_r$ is
a transitive, reflexive partial ordering, extending $\leq_T$, Turing
reducibility, so that if $x \leq_r y$ and
$z \leq_r y$ then $x \oplus z \leq_r y$ holds, where 
 $x \oplus z$ is the recursive union of $x$ with $y$) we shall
use $\leq_r$ to denote also the partial ordering amongst $r$-degrees
so induced.

If $D$ is a set of $r$-degrees, we
 say that $\mb{\bm$e$}$ is an {\em $r$-minimal upper bound} for $D$ if
$$
\forall \mb{\bm$d$} \in D (\mb{\bm$d$} \leq_r \mb{\bm$e$}) 
\wedge \forall \mb{\bm$g$} \left( \forall \mb{\bm$d$} 
\in D( \mb{\bm$d$} \leq_r \mb{\bm$g$} \wedge \mb{\bm$g$} \leq_r 
\mb{\bm$f$}) \rightarrow \mb{\bm$f$} =_r \mb{\bm$g$} \right)
$$

By the expression ``$r^{\#}$ exists'' we mean that there is a closed
and unbounded class of indescernibles, $C^r$, for $L[r]$, the
constructible closure of $L$ with $r$, and by $O^{\dagger}$, the
existence of such a class of indescernibles for $L[\mu]$- an inner
model for a measurable cardinal -- if such exists. 
 By ``$\neg O^{\dagger}$'' we mean that such a class does not exist.
For information on games see \cite{Ka}.
By contrast with Theorem 1.1 Kechris had earlier shown \cite{K1}:
\thm \  $(ZFC + \mb{\bm$\D$}^1_2$-Determinacy$)$.  If $D
 = \{\mb{\bm$d$}_n\}_{new}$ is a countable collection of $\mb{$\D$}^1_3$
degrees, that all lie within a single $Q$-degree, then $D$ has 
a minimal upper bound.
\ethm

At first sight it is tempting to  conjecture:

{\em Conjecture 1}\ \  ($ZFC + \forall r$ ($r^{\#}$ exists))  Every
countable collection of $\D^1_3$-degrees has a minimal upper bound.

Whilst Theorem 1.2 lends weight to this assuming
$\mb{\bm$\D$}^1_2$-determinacy, there was a large gap between this
assumption and that of Thm 1.1.
 It is the purpose of this note to close  that gap.
We prove:

\thm \  Assume $ZFC + \forall r\;(r^{\#}$ exists$)$ . 
For every countable set $D = \{{\bm d_n}\}$ 
of $\D^1_3$-degrees, either $D$ has a
$\D^1_3$-minimal upper bound, or 
for every $n$, for any $r \in{\bm  d_n}$ $\D^1_2(r)$ Determinacy holds.
\ethm

The inner model machinery we use is due to Steel \cite{St1} building on
the fine-structure of iterations trees, due to Mitchell \& Steel
\cite{MiSt}.  We assume the reader is familiar with \cite{MiSt},
\cite{St1}, and \cite{W1}.
 The basic structure of the proof of Theorem 3 is still that of
\cite{W1}.  This note details how \cite{W1} may be amended or ``read''
to show how to relax  the extra  assumption of
\cite{W1}, and still make the proof of Theorem 1.3 go through.  We shall
occasionally 
try to avoid
wholescale repetition by assuming the reader has a copy of \cite{W1} to
refer to.

We first prove  fairly directly the following weakening of 1.3.
\thm Assume there are two measurable cardinals, and that
$\D^1_2$-Determinacy fails. Then every countable set of $\dew$
degrees has a minimal upper bound.
\ethm

We indicate at the end
how this assumption can be weakened to that of sharps for reals to get
1.3 

 \section{}

The basic tools of Theorem 1.1 were the $\Sigma^1_3$-correctness of the
Dodd-Jensen Core Model, $K_{DJ}$, together with the analysis of
$\Sigma^1_3$ sets as unions of $\aleph_1$-Borel sets of \cite{W2}, where
the codes of the Borel sets could be taken from reals coding wellorders
and mice, from $K_{DJ}$.  Essentially this analysis held because the
class of uniform indiscernibles $C $ for reals could be
computed in $K_{DJ}$ as
$$\stackrel{\disp{\bigcap C^r}}{r\subseteq \gamma < \omega^V_1}$$
 from sharps for bounded subsets of $\omega_1^V$ in $K$.  Thus ``$u_n =
u_n^{K_{DJ}}$'' for $n < \omega$.

Here we shall use:
\thm  (Steel)\cite{St1} 7.9   Assume there are two measurable cardinals $\kappa <
\Omega$, but that there is no inner model of a Woodin cardinal.  Then
the core model $K$ of \cite{St1} is $\Sigma^1_3$-correct. 
\ethm
We say that a model $M$ is $\S^1_3$-correct, if for any $x \in M$, and
any non-empty $\P^1_2(x)$ set $B$ of reals, then $B$ has a member $y \in M$.
It is still unknown whether $u_2 < \omega_2$ follows from the
hypothesis of Theorem 2.1, but its proof does show that
computing the uniform indiscernibles for bounded subsets of the lower
measurable cardinal $\kappa$, yields the same class, both in $V$ and
the Steel $K$.

Let $D = \{\mb{\bm$d$}_n\}_{n\in \o}$
 be a fixed countable collection of
$\D^1_3$-degrees.  Without loss of generality we assume $n < m
\rightarrow \mb{\bm$d$}_n <_3 \mb{\bm$d$}_m$.  Let $d_n \in
\mb{\bm$d$}_n \cap \baire$ 
be a representative of each class, and let $z = \oplus
\langle d_n | n < \omega \rangle$ some recursive coding of the $d_n$.

We assume for some $n$ that lightface $\D^1_2(d_n)$-Determinacy fails.
Clearly we may assume that $n = 0$, and in fact we shall take $d_0 =  0$
as  the proof
will relativize uniformly. 
Initially we shall work under the assumptions of Thm.  2.1, and at
the end of the paper indicate that they may be removed, if
we continue to assume the existence of $\#$'s for reals.
For $z$ a set of ordinals, we let $K^z$ be the Steel core model
relativized to the predicate $z$, and $(K^z)^c$ the associated
preliminary
model (see \cite{St1} \sec\sec 1 - 5).  All references to ``mice"
will mean an $\o$-sound, $\o$-iterable mouse 
as defined in \cite{MiSt}, Defs. 2.8.3 \& 5.1.4 (unless otherwise stated).

As $\D^1_2$-determinacy fails $K^c \models $ ``there are no Woodin cardinals"
({\em cf.} \cite{St1} 7.14) and 
all iteration trees then are {\em simple} ({\em cf.} \cite{MiSt} 5.1.2), 
and furthermore we have a
$\Pi^1_2$ iterability condition on premice.
\df  Let $f \in \baire$.  Let $\tau$ be such that
$$
\langle K^f_{\tau}, \in, f, E^{K^f} \rangle \prec_{\Sigma_1} \la
K^f_{\omega_1}, \in, f, E^{K^f} \rangle.
$$
Then $\tau$ is $K^f$-stable.  We set $\tau (f)$ to be the least $K^f$-stable.
\edf
\df  $\delta^1_3(f) = \sup \{ \| X\| : X \in WO \wedge X \leq_3 f\}$.
\edf
\lem If $K^f$ is $\S^1_3$-correct, $\tau (f) = \delta^1_3(f)$ for any $f \subseteq \omega$.
\elem
\pf Just as for $\delta^1_2$, the first stable ordinal and $L$, using the $\Sigma^1_3(f)$-correctness of $K^f$.
\q

In \cite{W1} we used mainly the fact that below a measurable cardinal
the constructible closure $K_{\omega_1}[f]$ had the same domain as
$K^f_{\omega_1}$.  We see no {\em a priori} reason for this to be true
for larger core models (we replace it with Lemma 2.18(iii) below). \mar However:
\lem If $\lambda \leq \omega^K_1$ and $y \in \baire$, then
$K_{ \lambda} =_{df}
\langle J^{E^K}_{{\bf \lambda}},
\in , E^K\rest \l \rangle = (K _ \lambda)^{K^y}$.
\elem
\pf Let $\overline{K} = K^{K^y}$.\\
{\em Claim}\ \ $\overline{K}$ is a universal weasel.   \\
\pf By the Weak Covering Lemma for $K^y$ (see \cite{MSS}) for a cub
class $D \subseteq \Omega$ of cardinals, we shall have 
$$ \beta \in D
\rightarrow cf (\beta^{+K^y}) \geq \beta. $$
 Similarily building $K$
inside $K^y$, (note that enough of the measure on $\Omega$
survives into $K^y$ for the construction to take place), 
we appeal to the Weak Covering Lemma for $K$ inside $K^y$
to get the same conclusion on a cub set $D_1 \subseteq \Omega$ of
$K^y$-cardinals

$$ \beta \in D_1 \rightarrow cf^{K^y}( \beta^{+{\overline{K}}}) \geq \beta.
$$

Consider the comparison of $K$ with $\overline{K}$ with resultant trees
${\cal T}, {\cal U}$. (Again we follow \cite{MiSt} \sec 7 in our
definition of comparison, excepting that we shall always
consider comparisons of mice as terminating in a common model).
 If $\overline{K}$ is not universal then, by a
standard argument,(\cite{St1}),\sec 8
 there is a cub  $E  \subseteq \Omega$, with $i < j
\in E \rightarrow \pi_{ij}^{{\cal T}}\  (\k_i) = \kappa_j = j$\ \ (where
$\kappa_i = crit(E_i^{{\cal U}}))$ with $\pi_{ij}^{{\cal U}}\ ''
\kappa_j \subseteq \kappa_j$,  whilst $\pi_{ij}^{{\cal T}}\ ''
(\kappa_i^+)^{M_i^{{\cal T}}}$ is cofinal in $(\kappa_j^+)^{M_j^{{\cal
T}}}$.  Letting $i_0 = \min E$, and choosing a regular
$j \in E \cap D \cap D_1$, we have $\p^\calu_{0,j}(\k_j) = \k_j$
and $cf(\p^\calu_{0,i}((\k^+_j)^{\bar K}))\geq \k_j$ for all $i \leq j$.
But the latter is $(\k_j^+)^{M^\calu_j}$
 and this has cofinality $ (\k_{i_0}^+)^{M^\calt_{i_0}} < \k_j$!
\q{(Claim)}

Consequently on neither side of the coiteration is there any truncation
on the main branch.  This can only mean $K _ \lambda  = \overline{K} _
\lambda$ for any $\lambda$ less than the first measurable of $K$ and of
$\overline{K}$. \q 

\df
For $x \in \mbox{}^\omega\omega\  \D^1_3(x) =
\{M\;|\; \mbox{$M$ a mouse, } M \in  K^x_{\delta^1_3(x)}\}$.
\edf

We shall  use a relativized form of the above argument:
\lem (i) For $x \in K^y \; x,y \in \mbox{}^\omega\omega\ \D^1_3(x) =
\left(\D^1_3(x)\right)^{K^y}$.
(ii) If $x \leq_3 y$ then $\D^1_3(x) \sset \D^1_3(y)$
\elem
\pf(i):
2.5, \mar relativized to $K^x$, shows that, if $\t = 
\o^{K^x}_1$,
then $K^x_\t = (K^x_\t)^{K^y}$. For (ii): by the same argument, noting 
that $x \in K^y$ and that $\d^1_3(x) \leq \d^1_3(y)$ (both by the
$\S^1_3(y)$-correctness of $K^y$). The result is immediate.
\hspace*{\fill\q}
\cor  If $x =_3 y$ then $\Delta^1_3(x) = \Delta^1_3(y)$.
\ecor

\lem Let $\g = \o_1^{K^f}$, 
then $K^f_{\g}\models \all x\ex \a [K^f_\a\models \ex M\ex\t(M \mbox{
a mouse and } x \in L_\t[M,f])]$
\elem
\pf In $K^f$ define $\bar K[f]$ (where $\bar K = K^{K^f}$). Let $\wt K^f
=(K^f)^{\bar K[f]}$. By the same argument as Lemma 2.5 \mar 
$|K^f_{\o_1}| = |\wt K^f_{\o_1}|\sset |\bar K[f]_{\o_1}| = |K^f_{\o_1}|$
and so the result is immediate.\q
Let $\delta = \sup_n \{ \delta^1_3(d_n) : n < \mb{\bm$\omega$}\}$, and 
$\d_n = \d^1_3(d_n)$.
We set $\D^1_3(D) = \cup_n \D^1_3(d_n)$.

 The following replaces the notion of ``$\lambda$-$K$-degree'' from
\cite{W1}.
\df\label{mdegree} Let $\lambda$ be a p.r. closed ordinal, $M$ mouse with $\l , On
\cap M$.  A $\l$-$M$-{\em degree} is an equivalence class of reals under
the relation ``$f\in L_\l[M,g]$'' which we write as $f \leq_{\l,M} g$. 
 \edf 
($L_\l[M,g]$ is the usual constructible closure of $M$ with $g$).
Clearly
$\leq_{\l,M}$ is a reducibility ordering.
\df Let $f \in \baire$. Let either $\D = \D^1_3(f), \bar \d = \d^1_3(f)$
or $\D = \D^1_3(D), \bar \d = \d$
be as defined above. Then
we set $x \leq _\D y$ if $\exists M \in \D $ with $x \leq_{\bar\d,M} y$.
\edf
It is not hard to see that
$\leq_{\D}$ is also a reducibility ordering (for transitivity one
needs to observe only that if $x \leq_{\D} y \an y \leq_{\D} z$
as witnessed by mice $M,N \in \D^1_3(f)$ (resp. $\D^1_3(D)$)
then there is a mouse $P \in
\D^1_3(f)$ (resp. $\D^1_3(D)$) constructibly coding them both by level
$\d$).
Similar remarks hold
for $\D = \dew(D)$. 
There are thus notions
of minimal $\l$-$M$, $\D^1_3(f)$, and $\D^1_3(D)$-degree {\em etc.}
Apart from the first though, they are of limited utility, as $\D^1_3(f)$
is too generalised a set of mice. We replace this notion with that
of $``\wdew(f)"$ defined below.

The proof of \cite{W1} used a forcing argument derived from Friedman
\cite{Fr} using perfect trees with various notions of pointedness. For any
of the above reducibility notions:
\df A perfect tree $T \subseteq 2^{<\omega}$ is $r$-{\em
pointed} if $ \forall f (f \in T \rightarrow T \leq_r f) $.
\edf
Here as elsewhere we shall write $f \in T$ to mean that $f$ codes
the characteristic function of a set of sequence numbers coding a branch
through $T$. Throughout this note $T$ and $T^*$ will refer
to perfect trees of sequence numbers.
The following is an entirely general fact about
pointedness:
\lem [Sacks \cite{Sa2} 2.3]\ \ a)\ If $T$ is $r$-pointed and $T^*
\subseteq T,\ T^*$ perfect, and $T^* \leq_r T$, then $T^*=_rT$, and
$T^*$ is $r$-pointed.\ \ b)\  If $T$ is $r$-pointed,
and $T\leq_r f$, then there is a perfect $T^* \subseteq T$,
$r$-pointed, and such that $T^* =_r f$. 
\elem
Let $M$ and $\l$ be as in  definition \ref{mdegree}. Part a) of the next lemma
is just an application of the last fact with $\leq_r$ as $\leq_{\l,M}$.
Part b) is just a variation on the Sacks minimal degree construction
performed over $L[M]$.

\lem \label{sacks} a) ({\em cf. } [Sacks
\cite{Sa2} 2.3])  Let $T$ be $\l$-$M$-pointed.  Let $f\in
\mbox{}^\omega\omega$.  Then there is an $\l,M$-pointed $T^* \subseteq T$,
such that $T^*=_{\l,M}(T,f)$\ \  
b)({\em cf.} [Sacks \cite{Sa1}, 1.4 \& 3.1]) Let $T$ be perfect, $M$ a mouse with
$On \cap M < \lambda$, and suppose $g: \lambda
\stackrel{(1-1)}{\longrightarrow} \omega, \lambda$ p.r. closed.  Then
there is a perfect subtree $T^*$ of $T$, such that $T^* \in 
L_{\l + \o}[M,g]$
and $\forall f \in
T^* ( (f,T)$ is $\l$-$M$-minimal over $T$). 
\elem

\df  A premouse $N$ is {\em reasonable} if there is an 
$N$-cardinal, $\l = \l_N > \o$, with $H_\l^N$ closed under the sharp operation.
\edf

Let $S \subseteq \omega \times ^{\omega}\omega$ be a universal
$\Sigma^{1}_{3}$ set; and further require that for any 
$h \in \baire, \; S^{h} = \{e|S(e,h)\}$ is a complete
$\Sigma^{1}_{3}(h)$ subset of $\omega$.

Suppose $S(e,h) \longleftrightarrow \exists gP(e,h,g)$ where $P$ is
$\Pi^{1}_{2}$. For $e \in \omega$ and for $N$ reasonable, let $S^{N,e}_\l$
be the Martin-Solovay tree on $\omega^{2} \times u^{N}_{\omega}$
defined in $N$ whose projection is contained in $P$. (For definiteness,
let us take the definition of this ``MS"-tree as that of $S_2$ in
\cite{K3} \sec 2.2. Here we have set
$u^{N}_{1} = \lambda_{N}$, and $\langle u^{N}_{i}|i \leq \omega
\rangle$ enumerates the first $\omega$ members of

$$C^{N} = \bigcap_{a\subseteq\gamma<\lambda_{N}} I^{a}\backslash
\lambda_{N}$$

where $I^{a}$ is the class of Silver indiscernibles for $L[a]$.
This sequence is definable in $N$, and we
 use these indiscernibles to construct the tree  $S^{N,e}_\l$.

\rem (1)
A straightforward L\"{o}wenheim-Skolem argument using the full
Martin-Solovay tree in $K$ for $P$ constructed on the first
$\omega$-uniform indiscernible above the lower measurable $\k =
u^{K}_{1}$, shows that if $S(e,h)$ then there is a countable, reasonable
$N$ with some $g \in \baire$ so that $(g,h) \in p[S^{N,e}_{\l_N}]$.\\
(2) If $N \in M$ or $N = M$ are reasonable, and $\exists g (h,g) \in
p[S^{N,e}_{\l_N}]$, and $\l_N \leq \l_M$
then $\exists g (h,g) \in p[S^{M,e}]$.\\ 
\pf Clearly $H^N_{\lambda_N}
\subseteq H^M_{\lambda_M}$; hence by ``stretching'' the functions $f
\in N,\ f: [u^N_{\omega}]^k - u^N_{\omega}$ (each defined by $L[a]$
terms for $a \in H^N_{\lambda_N}$ with $a^{\#}$ existing in $N$) to
functions $\overline{f} \in M\ \overline{f}:[u^M_{\omega}]^k
\rightarrow u^M_{\omega}$ we get the desired conclusion.

\df
(i) $F(e,h) =_{df}$\\
\indent$ \{N\;|\;N 
\mbox{ is a } \leq_*\mbox{-least reasonable mouse so that } \ex g S^{N,e}_{\l_N}
(h,g) \mbox{ is illfounded }\}$. 
\\(ii) $\widetilde \D^1_3(h) = \bigcup\{ F(e,\bar h) \cap \dew(h)
 \;|\; {e\in \o} \an \bar h \leq_3 h\}$.
\edf
(i) above fulfills the role of $F(e,h)$ of \cite{W1}. 
 Here $\leq_{*}$ is the natural mouse
ordering: $M <_{*} N$ if the comparison process, {\em via} $\o$-maximal
iteration 
trees ${\cal T}, {\cal U}$ to a common model $M^\calt_\th =
M^\calu_\th$, then there has been a `drop' in degree along the main
branch $[0,\th]_U$ or in model, {\em i.e.} $D^\calu \cap [0,\th]_U \neq \emp$.
The next lemma shows that the mice of
$\wt\D^1_3(h)$ are thus the ``witnessing mice" for the $\S^1_3(h)$ complete
set. For the Dodd-Jensen core model $\wt\D^1_3(h)  = \D^1_3(h) $ although
there seems no {\em a priori} reason for this to be true here.

\lem $N \in F(e,h) \equi \\
 L[N,h] \models ``
\ex g S^{N,e}_{\l_N}
(h,g) \mbox{ is illfounded } \an \all \a V^{Coll(\o,\a)}\models`` \all M <_* N
\all g  S^{M,e}_{\l_M}
(h,g) \mbox{ is wellfounded" }$".  $``N \in F(e,h)"$ is thus a $\S^1_3(N,h)$
relation.
\elem
\pf ($\leftarrow$) 
If $F(e,h) \neq \emp$ but $N \notin F(e,h)$, this is because there is $M <_* N$ with
$\ex g S^{M,e}_{\l_M}
(h,g)$  is illfounded. Let $(M,N)$ be compared with resulting trees
$\calu,\calt$. Then using Remark (2) above, we see, if $P = M^\calu_\th$, the
last model on $\calu$ with $D^\calu\cap[0,\th]_U = \emp$, 
that if $\g = i^\calu_{0,\th}(\l_M)$, 
$\ex g S^{P,e}_{\g} (h,g) $ is illfounded,
 we may
map up the branch, sequence by sequence,
 to one in $S^{P,e}_\g$ using the iteration map
$i^\calu_{0,\th}$.
Let $\a = On \cap M$, and then by Shoenfield absolutness, if $G$ is $Coll(\o,\a)$-
generic over $L[N,h]$ we have $L[N,h,G]\models ``\ex \mbox{ countable iteration tree
 $ \calu $ on some premouse $M$ with }  [0,\th] \cap D^\calu = \emp \,\an\, 
\ex g S^{M^\calu_\th,e}_{\l} (h,g) \mbox{ is illfounded}$." This contradicts 
our assumptions.
($\rightarrow$) is 
straightforward. \q

\lem (i) If $e \in S^h$ then $F(e,h)\cap \dew(h) \neq \emptyset$. (ii)
$\wdew(h)$
is $\leq_{*}$ cofinal in
$\D^{1}_{3}(h)$. \\
(iii) $\all M \in \D^1_3(h)\ex N \in \widetilde
\D^1_3(h)(M \in
L_{\d^1_3(h)}[N,h])$. Hence $|K^h_{\d^1_3(h)}| = \bigcup_{N \in  \wdew(h)}|L_
{\d^1_3(h)}[N,h]|$.\\
(iv)  $\all M \in \D^1_3(h)\ex N \in \widetilde
\D^1_3(h)\all N' \in \dew(h) (N' \geq_* N \rightarrow M \in
L_{\d^1_3(h)}[N,h])$.
\elem 
\pf
(i) follows from the last sentence of Lemma 2.17 \mar since
$e \in S^h \equi \ex \a K^h_\a\models \ex N \in F(e,h) \equi
\ex \a< \d^1_3(h) K^h_\a\models \ex N \in F(e,h) $.

Suppose for (ii) the given set is not $\leq_*$-cofinal but there is $M_0
\in \D^1_3(h)$ such that
$\forall e( e \in S^h \rightarrow
[\forall N \in F(e,h)  N \leq_* M_0]$.

{\em Claim}\ \ ``$e \in S^h$'' is a $\D^1_3(h)$ relation
(contradicting the assumption on $S^h$ that it 
is a complete $\Sigma^1_3(h)$ set).\\
\pf \\
(1) $e \in S^h \equi \ex \calt\in K^h_{\d^1_3(h)}$, a countable iteration
tree $\calt$ on $M_0$ with last model $M^\calt_\infty$ so that 
\indent ``$\ex g S^{M^\calt_\infty,e}_\l(h,g) \mbox{ is illfounded }$".\\
\pf Suppose $e\in S^h$ and let $N \in K^h_\d$ where $\d = \d^1_3(h)$,
witness this using (i) above. Let $\calu,\calt$ be the trees of length
$\th$ resulting
from the comparison of $N,M_0$. As the latter are both countable,
$\calu, \calt \in K^h_{\o _1}$, and hence in $K^h_\a$ for some $\a <
\d$.
As $N \leq_* M_0 \; D^{\calu}\cap[0,\th] = \emp$, and so,
 if $S^{N,e}_\l$
is the tree defined in $N$, and
$i^\calu_{0,\th}(S^{N,e}_\l) = S^{ P,e}_\g$ where $P = N^\calu_\th$ is the last
common model on $\calt$ we have $D^\calu \cap [0,\th] = \emp$. If
$\g = i^\calu_{0,\th}(\l)$, and if $(g,h,f) \in [S^{N,e}_\l]$, we may again
map up the branch, sequence by sequence,
 to one in $S^{P,e}_\g$ using the iteration map
$i^\calu_{0,\th}$. 
Hence $(g,h) \in p[S^{P,e}_\g]$, and as $M^\calt_\th = M^\calu_\th = P$
we have shown the right hand side.
The converse direction is immediate. \qd{(2)} 

By Shoenfield absoluteness,
we have that the following is true:\\
(2) $e \in S^h \equi$ \\
 $\ex \a \;(V^{Coll(\o,\a)})^{L[M_0,h]} \models
``\ex \mbox{ a countable iteration tree $\calt$ on $M_0$ so that }
L[M^\calt_\infty,h]
\models ``e \in S^h"\; "$ 

Hence ``$e \in S^h$" is computable in $L[M_0,h]$, and we may write $e
\in S^h \equi n_0(e) \in (M_0,h)^\#$ for some recursive $n_0 \in
\baire$. But the latter sharp is in $K^h_\d$ - as $M_0, h$ are. \qd{(ii)}

To see that (iii) holds, let $M \in \D^1_3(h)$. Let $X \in K^h_\d$ be a
$\D^1_3(h)$ code for $M$ with $X \sset \o$. Let $e_0,e_1 \in \o $ be such
that $ n \in X \equi R(e_0,n) \equi \neg  R(e_1,n)$ where $R $ is a
universal $\S^1_3(h)$ subset of $\o \times \o$. Let $\a < \d$ be such
that $X \in K^h$ and that $\ex \bar N \in K^h_\a \ex \l$ an $\bar N$-cardinal, 
 with $K^h_\a \models `` \la e_0,n\ra \in R  \equi L[\bar
N,h] \models ``  \la e_0,n\ra \in p[S_\l^{R}]"\; "$, where $S_\l^{R}$
is a version of the MS tree for the relation $R$ defined in $\bar N$. 
Let $N \geq_* \bar
N$ with $N \in F(e,h)$ for some $e \in \o$. But now argue as at (2) in
part (ii) above that \\
(3):
$ \la e_0,n\ra \in R \equi$\\
  $\ex \b(V^{Coll(\o,\b)})^{L[N,h]} \models
``\ex \mbox{ a countable iteration tree $\calt$ on $N$ so that }
L[M^\calt_\infty,h]
\models ``\la e_0,n\ra \in R"\; "$ \\
hence $X$ and so $M$ is computable in $L_{\d^1_3(h)}[N,h]$. For the last sentence of (iii)
just note that what we did for $M$ we could have done for any $Y \in
|K^h_{\d^1_3(h)}|$. 

For (iv) note that we made no intrinsic use of ``$N \in F(e,h)$" in (3):
any $N' \geq_* \bar N$ would do. Since $N' \in \dew(h)$ we can bound the
$``\ex \b"$ computation, as there, by $\d^1_3(h)$. \q

We define a reducibility ordering to replace that of $``\dew(T)"$ from
\cite{W1}:
\df Set $f \leq_{\wdewt} g \Equi\quad  \ex \vec M = M_0,\ldots,M_k \in \wdewt
\quad f \in L_{\d^1_3(T)}[\vec M,g]$
\edf
Again there are corresponding notions of $\wdewt$-degree, minimality,
and perfect $\wdewt$-pointed trees.

\lem ({\em cf.} \cite {W1} Lemma 9) Let $T$ be a perfect tree, $f \in  T$
and let $F(e,f) = N$ for some $e \in  \omega$, but $f \notin
K^T$.  Then there are a
perfect $T^* \subseteq T,\ T^* \leq_3T$, and $M \in  \D^1_3(T)$
such that $\forall f \in  T^* (F(e,f) \neq \emp \an \all P \in F(e,f)\,
P \leq_* M)$. 
If additionally $T$ is $\wdewt$-pointed and $\Delta^1_3$-pointed then also
 $T^* =_3T$ and $T^*$ is $\wdew(T^*)$-pointed and $\dew$-pointed.
\elem
\pf Let $N \in  F(e,f)$.   Let $\bar N \in K^T$
be $\leq_*$-least with $\ex g S^{\bar N,e}_\l(f,g)$ illfounded
 for some $\bar N$-cardinal $\l$. Then  $\bar N \geq_*N$. Let $\a = On
\cap N$. 
Let $G$ be a $Coll(\o,\a)$-generic collapse over $L[\la T,\bar N\ra
^\#]$,
and so over $L[\la \bar N,T\ra]$, 
 with $f \not\in  L[ G,\la T,\bar N\ra ^\#]$. From $G$ define $h \sset
\o\times \o$ coding $\bar N$, with $f \notin L[\la h,T\ra^\#]$.
Thus
\begin{eqnarray*}
(1)\mbox{  \hspace*{0.3in}}\ex h(h \mbox{ codes a mouse} N_h \an \ex f \in T( 
\ex g S^{N_h,e}_\l(f,g) \mbox{ illfounded } \an 
f \notin L[\la h,T\ra^\#] \\
\an \all \bar h (( \bar h \mbox{ codes a countable premouse } M \an M \in
L[\la h,T\ra] \models ``M <_* N") \rightarrow \all g S^{M,e}_\l(f,g)
\mbox{ is wellfounded }]])&&
\end{eqnarray*}

This is $\S^1_3(T)$ (noting that there is a set of codes $\bar h$ of
countable
premice satisfying the last conjunct which is a recursive in $\la h,T\ra^\#$
set, $W_\e^{\la h,T\ra^\#}$ say, for some index $\e \in \o$).
 Hence there
is such an $h,N_h \in K^T_{\d^1_3(T)}$, and such an $f$ 
as in (1) with $f \notin
L[\la h,T\ra^\#]$. For this $h,N_h$ then 
\begin{eqnarray*}
E = \left\{ f\;|\; f\in T \an \ex g S^{N_h,e}_\l(f,g) \mbox{ is illfounded }
\an \all \bar h (\bar h \in W_\e^{\la h,T\ra^\#} \rightarrow
\all  g S^{M_{\bar h},e}_\l(f,g) \mbox{ is wellfounded })\right\}
\end{eqnarray*}
Then $E$ is $\S^1_2(\la h,T\ra^\#$, and contains elements not in
$L[\la h,T\ra^\#]$.
  Hence there is a perfect set
(given by a $T^* \subseteq T$) of such $f$, recursive in $\langle
h,T\rangle^{\#\#}$.  Hence $T^* \leq_3 (h,T) \leq_3 T$.  
For the last part, if $T$ is $\dew$-pointed, $T \leq_3 T^*$.
$\wdew(T^*)$-pointedness of $T^*$ is then immediate (as is 
$\dew$-pointedness)\q

\lem (cf. $[F]$ Lemma 8.)\ \ \label{extension} Let $T$ be 
$\D^1_3$-pointed, 
and  $\wdewt$-pointed, and suppose $f \in \omega^{\omega}$.  Then there
is a $\wdewt$-pointed $T^* \subseteq T$ such that $T^*
=_3(T,f)$.  Any $T^*$ satisfying these conditions is
$\D^1_3$-pointed and $\wdew(T^*)$-pointed.
\elem
\pf
 By general pointedness arguments, since $T$
is itself $\wdew(T)$-pointed, there exists a $\wdew(T)$-pointed $T^* \sset T$ such that\\
(1) $T^*
=_{\widetilde\D^1_3(T)} (T,f)$.\\ 
 Consequently,  setting $\bar \d = \d^1_3(T)$, 
for some finite sequence of mice $\vec M \in \wdew(T)$,
$T^* \in L_{\bar\d}[\vec M,(T,f)] \sset K^{(T,f)}_{\d^1_3((T,f))}$, so
$T^* \leq_3(T,f)$. Conversely note that $T \leq_3 T^*$ (since $T \leq_3
g$ where $g$ is the leftmost branch of $T^*$ which is recursive in
$T^*$). As $\wdewt \sset\D^1_3(T)\sset\D^1_3(T^*)$ (1) yields
$(T,f) \leq_3 T^*$. As $\wdewt \sset \wdew(T^*)$ and $\bar \d \leq \d^1_3(T^*)$
we have $T^*$ is $\wdew(T^*)$-pointed.

 As  $g \in T^* \rightarrow g \in T$ and
$T$ is
$\D^1_3$-pointed, $T \leq_3 g$, and then
 $T \in K^g_{{\delta}^1_3(g)}$.  Thus
$\widetilde\D^1_3(T) \subseteq \D^1_3(g)$.  
As $T^*$ is $\widetilde\D^1_3(T)$-pointed, for some mice $\vec M \in \dew(T)$
 $T^* \in L_\d[\vec M, g]$. But the latter is contained in
$ |K^g_{\delta^1_3(g)}|$, by the above.
Hence $T^* \leq_3 g$ and so is $\dew$-pointed.   \q

\lem Suppose additionally in Lemma \ref{sacks} b), that 
$T$ is $\wdewt$-pointed and $\D^1_3$-pointed,
and that $\omega < \lambda <
\delta^1_3(T)$ is p.r. closed.  Then, if $M \in \D^1_3(T)$ with $On \cap M
< \l$, the
$T^*$ of the conclusion of b) can be taken 
 to be $\widetilde \D^1_3(T^*)$ pointed and 
with $T^* =_3 T$ (and so $\dew$-pointed).
\elem
\pf
Set $\g = \delta^1_3(T)$.  
$K_\g^T \models \forall
\alpha \ \overline{\overline{\alpha}} = \omega$.  As $\lambda <
\g$ pick $g\in K_{\g}^T, g: \lambda 
\stackrel{(1-1)}{\longrightarrow} \omega$.
  Find $T^* \in L_{\l + \o}[M,g]$
as
in \ref{sacks} b) with the required $\l$-$M$-minimality.  
Then $T^* \in K^T_{\d^1_3(T)}$, so $T^* \leq_3 T$. As $T^*\sset T$, $T^*$
is $\dew$-pointed. As $f$, the leftmost branch of $T^*$ is recursive in $T^*$, we have $T\leq_3 T^*$,
so $T =_3 T^*. $
As $T$ is $\wdewt$-pointed, if $f \in T^*$, then  $T \in L_\g[\vec M,f]$ for
 some $\vec M \in \wdewt$. As $M,g \in K^T_\g$, by \mar 2.18 (iii), there
is $\vec N \in \wdewt$ with $(M,g)\in L_\g[\vec N,T]$. As $T^* \in L_\g[M,g]$
we have $T^* \in L_\g[\vec M, \vec N, f]$ and so $T^*$ is $\wdewt =\wdew(T^*)$-
pointed. \q

The following lemmas are used to control the growth of $\dew(f)$ in our
construction.
They are used
just as \cite{Fr} Lemmas 2 \& 3 are. 

\lem \label{A} Suppose a) $f$ is $\dew(D)$-minimal over $D$, b) $\dew(D)
\sset\dew(f)$, and c) $\all M \in\wdew(f)\ex N \in \dew(D) \;(M \leq_* N)
$. Then $f$ is $\dew$-minimal over $D$.
\elem
\pf By b) clearly $\all n (d_n \leq_3 f)$. And by a) $\all n (f \not \leq_3
d_n)$ (for suppose $f \in M = K^{d_n}_\a$ for some $\a < \d^1_3(d_n)$; then
$f$ is not $\dew(D)$-minimal over $D$!) Now suppose $\all n (d_n <_3 h)
\an h \leq_3 f$. Then $\dew(D) \sset \dew(h) \sset \dew(f)$. Then,
setting $\D = \dew(D)$: \\
(1) $\all n (d_n \leq_\D h)$.\\
But $h \leq_\D d_n \rightarrow \ex k \ex M \in \dew(d_k)
\ex \a < \d^1_3(d_k)\; h \in L_\a[M,d_n]$. Letting $m = \max\{n,k\}\;
h \in L_{\d^1_3(d_m)}[M,d_m]$. Thus $h \leq_3 d_m$, contradicting our 
supposition. Hence:\\
(2) $\all n ( d_n <_\D h)$.\\
(3) $h \leq _\D f$.\\
$ h \in K^f_{\d^1_3(f)}$ implies by Lemma 2.18(iii) \mar 
$\ex M \in \widetilde \dew(f)\ex \a < \d^1_3(f)\; h \in L_\a[M,f]$. By
b), c) and 2.18(iv) we can pick $N \geq_* M, \, N \in \D \sset \dew(f)
\an h \in L_{\d^1_3(f)}[n,f]$.
Hence (3) holds.
By a) then, $f \leq_\D h$. So $f \in L_\a[M,h]$ some $M \in \D \sset \dew(h)$
and $\a < \d^1_3(h)$. Hence
$f \in   
K^h_{\d^1_3(h)}$. Thus $f \leq_3 h$ as required. \q

For the next lemma as we have $\dew(d_n) \sset \dew(d_{n+1}) \sset \dew(D)$,
we assume we have 
 $\la M_n\;|\; n < \o \ra$ an  enumeration of $\dew(D)$ in such
a way that $\all i,j < \o \ex m > i,j \; M_i,M_j \in L_\l[M_m]$ where
$\l$
is a p.r. closed  ordinal $ > On \cap M_m$. So let $\la \l_m \;|\; m < \o\ra$
be an ascending sequence of p.r. closed  ordinals  with $\l_m > On \cap
M_m$, so that  $\l_m$ witnesses this, 
 and with $\sup_m\l_m = \d$. By thinning out the $\la {\bm d_n}
\;|\; n < \o \ra $ sequence, we may assume that for any $n$ $\l_n <
\d^1_3(d_n)$. 

\lem \label{B} Suppose a) $\all k \; (f,d_k)$ is $\l_k$-$M_k$-minimal over
$d_k$; b) $ \dew(D)\sset\dew(f)$; c) $ \dew(D)\supseteq\wt\dew(f)$; and 
d) $\all n \;(d_n <_3 f)$. Then $f$ is $\dew$-minimal over $D$.
\elem
\pf Set $\D = \dew(D)$. First note that $\all n(d_n <_\D f)$ (by b) \& d),
and using the argument of (2) of Lemma \ref{A}). We show that
$f$ is $\D$-minimal, and then the result follows by \ref{A}. So suppose\\
e) $\all n (d_n <_\D h)$ and \\
f) $h \leq _\D f$.\\
By f), for some $k$\\
(1) $h \in L_{\l_k}[M_k,f]$ (using our presumed properties on $M_m$ \& $\l_m$).\\
Similarly e) shows:
(2) $d_k \in L_{\l_n}[M_n,h]$ for some $k \leq n < \o$.\\
Now consider $(h,d_k)$. $d_k  \in L_{\l_k}[M_k,(h,d_k)]$ 
whilst $(h,d_k) \not  \in L_{\l_k}[M_k,d_k]$ by e). By (1) 
$(h,d_k)   \in L_{\l_k}[M_k,(f,d_k)]$. By $\l_k$-$M_k$-minimality of
$(f,d_k)$ over $d_k$, we conclude 
$(f,d_k)   \in L_{\l_k}[M_k,(h,d_k)]$. Using (2) we have 
$f \in L_{\l_n + \l_k}[M_n,M_k,h]$. For some sufficiently large 
$m > n,k,\; f \in L_{\l_m}[M_m,h]$. That is, $f \leq_\D h$. Hence
$f$ is $\D$-minimal as required. \q

We now have set up all the machinery, to
run the main argument of \cite{W1} Lemma 13,
 keeping roughly to the same notations.  

We state this as follows:
\lem Let $D, \, \{d_n\},\, \d, z, \l_n, M_n$ be as above. There is a perfect
set $T_0 \in K^z$ so that a) $\all f \in T_0\all n(d_n \leq_3 f)$;
b) $ \all f \in T_0 (f \notin \bigcup_nK^{d_n} \rightarrow 
f \mbox{ is } \dew\mbox{-minimal over } D \an \all M \in \wdew(f)\ex N
\in  \dew(D)\, M \leq_*N)$.
\elem\pf We only sketch the construction which takes place 
in $K^z$, defining  a binary system
of $\dew$-pointed trees $\la T_s | s \in 2^{<\o}\ra$.  Let $e_i$
be an enumeration of $\o$ in which every integer ocurs infinitely often.
Let $T_{\norwegian} = 2^{<\o}$. Assume $T_s$ has been defined $\all s \in 2^{<\o}\,
(lh(s)\leq i)$ so that \\
\indent
a) $lh(s) = j \leq i \rightarrow T_s =_3 d_j$\\
\indent
b) $T_s$ is $\wdew(T_s)$-pointed and $\dew$-pointed.\\
For $lh(s) = i$, define disjoint $T_{s\smallfrown 0}, T_{s\smallfrown 1}$
disjoint subtrees of $T_s$ according to the following recipe:\\
1) Split $T_s$ into two disjoint subtrees $T^*_s =_3 T^{**}_s =_3 d_{i+1}$,
perforce both $\wdew(T_s)$-pointed,\\
2) Then find $T^o_s \sset T^*_s, \; T^o_s\, \wdew(T^o_s)\mbox{-pointed }$, with
$T^o_s =_3 T^*_s =_3 d_{i+1}$ and so that 
$$ \all f \in T^o_s((f,T^o_s) \mbox{ is $\l_{i+1}$-$M_{i+1}$-minimal over }
 T^*_s)
$$
Define $T^{oo}_s$ entirely similarly using $T^{**}_s$, replacing $\mbox{}^*$
with $\mbox{}^{**}$ and $\mbox{}^o$ with $\mbox{}^{oo}$ throughout.\\
3) Find $T_{s\smallfrown 0} =_3 T^o_s =_3 d_{i+1} =_3 T^{oo}_s =_3 
T_{s\smallfrown 1}$, $T_{s\smallfrown j} \, \wdew(T_{s\smallfrown j})$-pointed, so that 
for some $N \in \widetilde \dew(T_s)\; \all f \in T_{s\smallfrown 0}
(\all M \in F(e_i,f), M \leq_*  N)$ if such a $T_{s\smallfrown 0}$ exists. Otherwise set
$T_{s\smallfrown 0} = T^0_s$. Similarly define $T_{s\smallfrown 1}$.
This induction can take place in any  initial segment of $K^z$, 
which is a model
of  $ZFC^-$, containing $z, \{K^{d_n}_{\d^1_3(d_n)}\}$
 (and so the $\wdewt(d_n)$
etc.)

The perfect 
tree $T_0$ of the lemma is that arising from the fusion of the
 $\la T_s | s \in 2^{<\o}\ra$, and so $[T_0] = \{f\;|\;
\bigcup_{g\in 2^\o}\bigcap_{N\in \o} f\in [T_{g\rest 
n}\}$.  The argument that $T_0$ is properly 
defined follows that of \cite{W1} with minor changes (replacing
$K[d_n]$ with $K^{d_n}$ for example, and replacing
$\dew(T)$ by $\wdewt$ throughout). 
By asking for $\wdew(T_s)$
pointed trees at each stage we have ensured
 the construction is absolute between 
$K^z$ and $V$, and one may show by induction
 on $lh(s)$,  using 2.20-2.22, that
 that the trees are 
in fact $\dew$-pointed in $V$.
We have ensured 
$ \all f \in T_0\all i(lh(s) = i \rightarrow
((f,T^o_s) \mbox{ is $\l_{i+1}$-$M_{i+1}$-minimal over }
 d_{i+1}$ (a notion absolute between $K^z$
and $V$) to fulfill requirement a) of Lemma \ref{B}. As clearly $d_n
\leq_3 f$ for any $f \in T_0$ (by $\dew$-pointedness of $T_s$ where
$f \in T_s \an lh(s) = n$), we have requirement b) of \ref{B} (and a)
of the current lemma). Now if $f \in T_0 \an f \notin \bigcup_n
K^{d_n}$,
by Lemma  2.20, if $F(e_i,f) \neq \emp$
then $\all M \in F(e_i,f)\, M \leq_* N$ for some $N \in \dew(D))$.
We have that $f,D$ satisfy 2.23 c) and the result follows by that lemma. \q
The following theorem analogous to \cite{W1} Thm 14 is proven
similarly {\em mutatis mutandis}.
\thm (i) For every $f \in T_0$, $f$ is an upper bound for $\{{\bm d_n}\}_{n
\in \o}$; if additionally $f \in T_0\backslash \bigcup_{i\in \o}K^{d_i}$
then ${\bm f}$ is a minimal upper bound for $\{{\bm d_n}\}_{n
\in \o}$.\\
(ii) There exists a minimal upper bound of $\{{\bm d_n}\}_{n
\in \o}$.\\
(iii) There is a least upper bound of $\{{\bm d_n}\}_{n
\in \o}$ iff $\ex i_0 \in \o(\re \sset K^{d_{i_0}})$.
\ethm
We remark now on how to remove the assumption that there are two
measurable cardinals $\k < \O$. 
The essential ingredient of Lemma 2.25 \mar is to have a $ZFC^-$ model
$M$, containing $z$ which is $\S^1_3$-correct (and so contains all
$\dew(d_n)$ - thinking of this as the mice coded by $\dew(d_n)$ subsets
of $\o$), and so has some ordinal height $\th > \sup_n\{\d^1_3(d_n)\} =
\d$). Given then an enumeration of $\dew(D) = \bigcup_n \dew(d_n)$, and a
sequence $\l_n < \d^1_3(d_n)$ with $\sup \l_n = \d$, satisfying the
requirement before Lemma 2.24 \mar (which we can define over $M$), we
define a sequence of trees $\la T_s\ra$ as above by induction in $M$. We
establish the existence of such a suitable $M$ by the following form of
argument (due to Woodin? - this is the kind of argument that Hauser uses
to lift the theorem  of \cite{W2}, on all $\P^1_3$ sets of reals containing 
$\P^1_3$-singletons being equivalent with all reals residing in the
Dodd-Jensen $K_{DJ}$ to the context of the Steel Core Model
- we should like to thank him for explaining this argument to us.)

We may assume $\all y\, y^\dagger$ exists (otherwise for some $ y \geq_3
z$
satisfying $\neg y^\dagger$ we could run the argument above using the
model $K_{DJ}[y]$ as is done in \cite{W1}.) But then, we can consider a
(or any) canonical inner model arising from $y^\dagger$ with a measure
$\mu$ on some $\k > \o_1$ containing $y,\; L[\mu,y]$ say, and using as $\O$
any of the upper indiscernibles from $y^\dagger$ and as in \cite{St1}
7.7, construct a $\S^1_3$-correct model $P = K^z$ in $L[\mu,y]$. Let
$\t = \o_1^P$, and let $M(y) = P|\t$. Note that for any $y' \geq_T y$ a
trivial comparison argument shows $M(y)$ an initial segment of $M(y')$.
Let $M$ be the union over all such $y \geq _T z$. The following
claim shows that $M$ has sufficient properties for the induction to
go through.\\
{\em Claim}  (i) $ z \in M \models ZFC^-$; (ii) $M$ is $\S^1_3$-correct.\\
For (i), clearly $M$ is of the form $\la J^{E^z}_\th,\in, E^z\ra$ for
some $\th \leq \o_1$.  So assume $\th < \o_1$ and that $ZFC^-$ fails.
Let $w$ code $M$ and place ourselves in a $L[\mu,w]$. Again a trivial
comparison shows us that $M$ is a proper initial segment of $K^z_{\o^{K^z}_1}$,
but $M$ is the union of such! For (ii), if $t \in J^{E^z}_\a$ (the
latter
coded by some real $y$ say) and if $B$
is a non-empty set of $\P^1_2(t)$ reals, let $s \in B$, and place
ourselves
in $L[\m,s\oplus y]$: $t \in J^{E^z}_\a$ is an initial segment of $K^z$
and the latter is $\S^1_3$-correct. Hence there must be an $s' \in B$
with
$s' \in K^z$. Hence $s' \in M$ as required.


 \end{document}